# A NOTE ABOUT KHOSHNEVISAN–XIAO CONJECTURE

By Martynas Manstavičius

*University of Connecticut*

Khoshnevisan and Xiao showed in [*Ann. Probab.* **33** (2005) 841–878] that the statement about almost surely vanishing Bessel–Riesz capacity of the image of a Borel set $G \subset \mathbb{R}_+$ under a symmetric Lévy process $X$ in $\mathbb{R}^d$ is equivalent to the vanishing of a deterministic $f$-capacity for a particular function $f$ defined in terms of the characteristic exponent of $X$. The authors conjectured that a similar statement is true for all Lévy processes in $\mathbb{R}^d$. We show that the conjecture is true provided we extend the definition of $f$ and require certain integrability conditions which cannot be avoided in general.

**1. Introduction.** Given a Borel set $G \subset \mathbb{R}^d$, let $\mathcal{P}(G)$ be the set of all Borel regular probability measures on $G$. Let $f:\mathbb{R}^d \to [0,\infty]$ be a Borel measurable function. Then for any $\mu \in \mathcal{P}(G)$, one can define the $f$-energy for this measure $\mu$ as

$$\mathcal{E}_f(\mu) := \int \int_{G \times G} f(x-y)\, \mu(dx)\, \mu(dy).$$

Furthermore, the $f$-capacity of $G$ is defined as

$$\mathcal{C}_f(G) := \left[\inf_{\mu \in \mathcal{P}(G)} \mathcal{E}_f(\mu)\right]^{-1},$$

where we set $1/\infty = 0$. If the function $f(\cdot) = \|\cdot\|^{-\beta}$ for some $\beta > 0$, then $\mathcal{C}_f$ will be simply denoted by $\mathcal{C}_\beta$.

Now consider a Lévy process $X(t), t \geq 0$, in $\mathbb{R}^d$ with the Lévy triplet $(\mathbf{b}, A, L)$, where $\mathbf{b} \in \mathbb{R}^d$, $A$ is a $d \times d$ symmetric nonnegative definite matrix and $L$ is the Lévy measure, that is, a Borel measure on $\mathbb{R}^d \setminus \{0\}$ such that

$$\int_{\mathbb{R}^d \setminus \{0\}} 1 \wedge \|x\|^2 L(dx) < +\infty.$$

Then the Lévy–Khinchine formula is simply $E\exp\{i\langle\boldsymbol{\xi}, X(t)\rangle\} = \exp\{-t\Psi(\boldsymbol{\xi})\}$, where $\boldsymbol{\xi}\in\mathbb{R}^d$, $t\geq 0$, $\langle\cdot,\cdot\rangle$ denotes the usual scalar product of $\mathbb{R}^d$ and the Lévy exponent $\Psi(\cdot)$ is given by

$$\Psi(\boldsymbol{\xi}) = -i\langle\mathbf{b},\boldsymbol{\xi}\rangle + \tfrac{1}{2}\langle A\boldsymbol{\xi},\boldsymbol{\xi}\rangle$$
$$- \int_{\mathbb{R}^d\setminus\{0\}} (e^{i\langle\mathbf{y},\boldsymbol{\xi}\rangle} - 1 - i\langle\mathbf{y},\boldsymbol{\xi}\rangle\mathbf{1}_{\{\|y\|\leq 1\}}(y))\, L(dy). \tag{1.1}$$

By [1], Theorem 2.2, we get that, for all $\beta\in(0,d)$ and all Borel sets $G\subset\mathbb{R}_+$,

$$\mathcal{C}_\beta(X(G)) = 0 \quad \text{a.s.} \tag{1.2}$$
$$\iff \quad \forall\mu\in\mathcal{P}(G): \int_{\mathbb{R}^d} \mathcal{E}_{\chi_{\boldsymbol{\xi}}}(\mu)\|\boldsymbol{\xi}\|^{\beta-d}\,d\boldsymbol{\xi} = +\infty,$$

where

$$\chi_{\boldsymbol{\xi}}(x) := e^{-|x|\Psi(\mathrm{sign}(x)\boldsymbol{\xi})} \qquad \forall\, x\in\mathbb{R}.$$

Note that this function $\chi_{\boldsymbol{\xi}}(x)$ is continuous both as a function of $\boldsymbol{\xi}$ and $x$, but its real part is not always positive and, therefore, not necessarily Lebesque integrable as a function of $\boldsymbol{\xi}\in\mathbb{R}^d$, as was kindly pointed out by the referee. In case the Lévy process $X$ is symmetric, one gets the following equivalency (see [1], Corollary 2.4):

$$\mathcal{C}_\beta(X(G)) = 0 \quad \text{a.s.} \iff \mathcal{C}_{f_{d-\beta}}(G) = 0, \tag{1.3}$$

where

$$f_\gamma(x) := \int_{\mathbb{R}^d} e^{-|x|\Psi(\boldsymbol{\xi})}\|\boldsymbol{\xi}\|^{-\gamma}\,d\boldsymbol{\xi} \qquad \forall\, x\in\mathbb{R}, \gamma\in(0,d). \tag{1.4}$$

Khoshnevisan and Xiao conjectured that (1.3) holds quite generally and we show this is indeed the case, subject to a slight change in the definition of the gauge function $f_\gamma(\cdot)$, as well as certain integrability assumptions that cannot be avoided in general. The details and appropriate examples are given in the next two sections.

**2. Results.** Given a general Lévy process on $\mathbb{R}^d$, we will need the real part of the function $\chi_{\boldsymbol{\xi}}(x)$ which, using the fact that the real part $\Re\Psi(\boldsymbol{\xi})$ is always even and the imaginary part $\Im\Psi(\boldsymbol{\xi})$ is always odd, can be written as

$$\Re\chi_{\boldsymbol{\xi}}(x) = e^{-|x|\Re\Psi(\boldsymbol{\xi})}\cos(|x|\Im\Psi(\boldsymbol{\xi})) = (\Re\chi_{\boldsymbol{\xi}}(x))_+ - (\Re\chi_{\boldsymbol{\xi}}(x))_-$$
$$= e^{-|x|\Re\Psi(\boldsymbol{\xi})}(\cos(|x|\Im\Psi(\boldsymbol{\xi})))_+ - e^{-|x|\Re\Psi(\boldsymbol{\xi})}(\cos(|x|\Im\Psi(\boldsymbol{\xi})))_-,$$

where we set $(h(x))_+ = \max\{h(x),0\}$ and $(h(x))_- = \max\{-h(x),0\}$. Now for any Borel set $G\subset[0,\infty)$ any $x\in G$ and $\gamma\in(0,d)$, define

$$g_{\gamma,+}(x) := \int_{\mathbb{R}^d} (\Re\chi_{\boldsymbol{\xi}}(x))_+ \|\boldsymbol{\xi}\|^{-\gamma}\,d\boldsymbol{\xi},$$
$$g_{\gamma,-}(x) := \int_{\mathbb{R}^d} (\Re\chi_{\boldsymbol{\xi}}(x))_- \|\boldsymbol{\xi}\|^{-\gamma}\,d\boldsymbol{\xi}. \tag{2.1}$$

Then we have the following:

THEOREM 2.1. *Let $X$ be a Lévy process in $\mathbb{R}^d$ with the Lévy exponent $\Psi$. Then for any Borel set $G \subset \mathbb{R}_+$ and any $\beta \in (0,d)$:*

(a) $\mathcal{C}_{g_{d-\beta},+}(G) > 0 \quad \Rightarrow \quad \mathcal{C}_{g_{d-\beta},-}(G) > 0.$
(b) $\mathcal{C}_{g_{d-\beta},-}(G) = 0 \quad \Rightarrow \quad \mathcal{C}_{g_{d-\beta},+}(G) = 0.$
(c) $\mathcal{C}_\beta(X(G)) = 0 \quad a.s. \quad \Rightarrow \quad \mathcal{C}_{g_{d-\beta},+}(G) = 0.$
(d) $\mathcal{C}_{g_{d-\beta},+}(G) = 0 \quad \not\Rightarrow \quad \mathcal{C}_\beta(X(G)) = 0 \quad a.s.$
(e) $\mathcal{C}_{g_{d-\beta},+}(G) = 0 \quad and \quad \int_{\mathbb{R}^d} \mathcal{E}_{(\Re\chi_\xi)_-}(\mu) \|\boldsymbol{\xi}\|^{\beta-d} d\boldsymbol{\xi} < +\infty, \quad \forall \mu \in \mathcal{P}(G) \quad \Rightarrow \quad \mathcal{C}_\beta(X(G)) = 0 \quad a.s.$

REMARK 2.2. If the Lévy process is symmetric, the Lévy exponent $\Psi$ is real-valued, and we always have $g_{\gamma,+}(\cdot) = f_\gamma(\cdot)$ and $g_{\gamma,-}(\cdot) \equiv 0$. Then (c) and (e) combined yield the result of [1], Corollary 2.4, since $\mathcal{E}_{(\Re\chi_\xi)_-}(\mu) \equiv 0$ for any $\mu \in \mathcal{P}$.

REMARK 2.3. The integrability condition of part (e) is rather stringent and it trivially holds for symmetric Lévy processes, so one might ask if there are nonsymmetric Lévy processes that satisfy it. This is indeed true as the example at the end of the paper shows.

**3. Proofs.** We begin with a simple observation:

LEMMA 3.1. *For any Lévy process $X$ in $\mathbb{R}^d$ and any $\mu \in \mathcal{P}(G)$, where $G \subset \mathbb{R}_+$ is a Borel set, we have*

(3.1) $$\mathcal{E}_{\chi_\xi}(\mu) = \mathcal{E}_{\Re\chi_\xi}(\mu).$$

PROOF. The proof is rather simple. As in the proof of Lemma 2.1 of [1], we get that, for any $s, t \geq 0$ and all $\boldsymbol{\xi} \in \mathbb{R}^d$,

$$\chi_{\boldsymbol{\xi}}(t-s) = E e^{i \langle \boldsymbol{\xi}, X(t) - X(s) \rangle}$$

and also, for any $\mu \in \mathcal{P}(G)$,

$$\mathcal{E}_{\chi_\xi}(\mu) = E \left| \int_G e^{i\langle \boldsymbol{\xi}, X(t) \rangle} \mu(dt) \right|^2 \in [0,1].$$

To see why (3.1) is true, simply write

$$\left| \int_G e^{i\langle \boldsymbol{\xi}, X(t) \rangle} \mu(dt) \right|^2 = \left( \int_G \cos(\langle \boldsymbol{\xi}, X(t) \rangle) \mu(dt) \right)^2 + \left( \int_G \sin(\langle \boldsymbol{\xi}, X(t) \rangle) \mu(dt) \right)^2$$

$$= \int\int_{G \times G} \cos(\langle \boldsymbol{\xi}, X(t) \rangle) \cos(\langle \boldsymbol{\xi}, X(s) \rangle) \mu(dt) \mu(ds)$$

$$+ \int\int_{G\times G} \sin(\langle \boldsymbol{\xi}, X(t)\rangle)\sin(\langle \boldsymbol{\xi}, X(t)\rangle)\mu(dt)\mu(ds)$$

$$= \int\int_{G\times G} \cos(\langle \boldsymbol{\xi}, X(t) - X(s)\rangle)\mu(dt)\mu(ds)$$

$$= \int\int_{G\times G} \Re e^{i\langle \boldsymbol{\xi}, X(t) - X(s)\rangle}\mu(dt)\mu(ds).$$

Taking expectations of both sides, using Fubini's theorem (which is justified since the integrand is bounded by an integrable function 1) and the fact that the operators $E$ and $\Re$ on the set of integrable complex-valued functions commute, we obtain (3.1). □

The next step is to establish the following:

LEMMA 3.2. *Let $g_{\gamma,\pm}(x)$ be as in (2.1). Then for any $\gamma \in (0,d)$,*

$$\int_{\mathbb{R}^d} \mathcal{E}_{(\Re\chi_{\boldsymbol{\xi}})_\pm}(\mu)\|\boldsymbol{\xi}\|^{-\gamma}\,d\boldsymbol{\xi} = \mathcal{E}_{g_{\gamma,\pm}}(\mu).$$

PROOF. Since both functions $(\Re\chi_{\boldsymbol{\xi}})_-\|\boldsymbol{\xi}\|^{-\gamma}$ and $(\Re\chi_{\boldsymbol{\xi}})_+\|\boldsymbol{\xi}\|^{-\gamma}$ are nonnegative, we can apply the Fubini–Tonelli theorem to get the result. □

We are now ready to establish the claims of the main theorem.

PROOF OF THEOREM 2.1. (a) and (c). Suppose $\mathcal{C}_{g_{d-\beta,+}}(G) > 0$. Then there exists a probability measure $\mu \in \mathcal{P}(G)$ such that $\mathcal{E}_{g_{d-\beta,+}}(\mu) < +\infty$. This and Lemma 3.2 yield

$$(3.2) \qquad \int_{\mathbb{R}^d} \mathcal{E}_{(\Re\chi_{\boldsymbol{\xi}})_+}(\mu)\|\boldsymbol{\xi}\|^{\beta-d}\,d\boldsymbol{\xi} < +\infty.$$

Recall that $\mathcal{E}_{\chi_{\boldsymbol{\xi}}}(\mu) \in [0,1]$ and, by Lemma 3.1, $\mathcal{E}_{\chi_{\boldsymbol{\xi}}}(\mu) = \mathcal{E}_{\Re\chi_{\boldsymbol{\xi}}}(\mu)$. Then using (3.2), we get both

$$\int_{\mathbb{R}^d} \mathcal{E}_{\Re\chi_{\boldsymbol{\xi}}}(\mu)\|\boldsymbol{\xi}\|^{\beta-d}\,d\boldsymbol{\xi} < +\infty \quad \text{and}$$

(3.3)

$$\int_{\mathbb{R}^d} \mathcal{E}_{(\Re\chi_{\boldsymbol{\xi}})_-}(\mu)\|\boldsymbol{\xi}\|^{\beta-d}\,d\boldsymbol{\xi} < +\infty,$$

since

$$\mathcal{E}_{(\Re\chi_{\boldsymbol{\xi}})_+}(\mu) = \mathcal{E}_{\Re\chi_{\boldsymbol{\xi}}}(\mu) + \mathcal{E}_{(\Re\chi_{\boldsymbol{\xi}})_-}(\mu).$$

The first integral in (3.3), together with (1.2), implies part (c) of the theorem, whereas the second proves part (a).

(b) This claim follows from part (a) and the fact that capacities are non-negative.

(d) All we need is an example. Consider $d=1$, $G=[0,1]$ and a deterministic process $X(t)=t$, $t\in[0,1]$. Then $\Psi(\xi)=-i\xi$ and $g_{1-\beta,\pm}(x)=+\infty$ for all $x\in\mathbb{R}\setminus\{0\}$, $g_{1-\beta,+}(0)=+\infty$, $g_{1-\beta,-}(0)=0$ and $\beta\in(0,1)$. Indeed, for any $\beta\in(0,1)$ and $x\neq 0$, we have

$$g_{1-\beta,+}(x) = \int_{\mathbb{R}}(\cos(|x|\xi))_+|\xi|^{\beta-1}\,d\xi$$

$$\geq \sum_{k=1}^{\infty}\int_{((4k-1)/2+1/6)\pi/|x|}^{((4k+1)/2-1/6)\pi/|x|}(-\cos(|x|\xi))\xi^{\beta-1}\,d\xi$$

$$\geq \frac{1}{2}\sum_{k=1}^{\infty}\int_{((4k-1)/2+1/6)\pi/|x|}^{((4k+1)/2-1/6)\pi/|x|}\xi^{\beta-1}\,d\xi$$

$$= \frac{\pi^{\beta}}{2\beta|x|^{\beta}}\sum_{k=1}^{\infty}\left[\left(\frac{4k+1}{2}-\frac{1}{6}\right)^{\beta}-\left(\frac{4k-1}{2}+\frac{1}{6}\right)^{\beta}\right]$$

$$= \frac{\pi^{\beta}}{2\beta|x|^{\beta}}\sum_{k=1}^{\infty}3^{-\beta}(6k-1)^{\beta}\left[\left(1+\frac{2}{6k-1}\right)^{\beta}-1\right]$$

$$\geq \frac{(2\pi)^{\beta}}{2|3x|^{\beta}}\sum_{k=1}^{\infty}(6k-1)^{\beta-1}=+\infty,$$

since the mean value theorem applied to the function $h(x)=(1+x)^{\beta}-1$, $x\in[0,1]$, yields

$$h(x)=h(x)-h(0)=h'(c)x\geq \min_{y\in[0,1]}h'(y)x=\beta 2^{\beta-1}x,$$

where $c\in(0,x)$. A similar argument gives $g_{1-\beta,-}(x)=+\infty$ for all $x\in\mathbb{R}\setminus\{0\}$. Values $g_{1-\beta,\pm}(0)$ are obvious. Now it is clear that $\mathcal{C}_{g_{1-\beta,+}}(G)=0$. Note also that we have shown $\mathcal{C}_{g_{1-\beta,-}}(G)=+\infty$ (one only needs to take $\mu=\delta_x$ for any $x\in G$, which leads to zero $g_{1-\beta,-}$-energy of a set $G$). At the same time, taking $\mu$ to be the Lebesgue measure on $G=[0,1]$, we have by direct computation

$$\mathcal{E}_{\beta}(\mu)=\int\int_{[0,1]^2}|t-s|^{-\beta}\,dt\,ds=\frac{2}{(1-\beta)(2-\beta)}<+\infty,$$

so $\mathcal{C}_{\beta}(X(G))=\mathcal{C}_{\beta}([0,1])>0$ for any $\beta\in(0,1)$.

(e) This is easy again, since the assumptions and the equality $\mathcal{E}_{(\Re\chi_{\xi})_+}(\mu)=\mathcal{E}_{\Re\chi_{\xi}}(\mu)+\mathcal{E}_{(\Re\chi_{\xi})_-}(\mu)$ yield

$$\int_{\mathbb{R}^d}\mathcal{E}_{\Re\chi_{\xi}}(\mu)\|\boldsymbol{\xi}\|^{\beta-d}\,d\boldsymbol{\xi}=+\infty \qquad \forall\mu\in\mathcal{P}(G).$$

Then Lemma 3.1 and (1.2) give the claim of (e).

This completes the proof of the theorem. □

EXAMPLE 1. To show that the conditions of part (e) of the theorem are satisfied not only for symmetric Lévy processes, consider a real-valued Poisson process $X_t$, $t \geq 0$, with parameter 1. The Lévy triplet [relative to the cut-off function $c(x) = \mathbf{1}_{\{\|x\| \leq 1\}}$ as in (1.1)] of such $X_t$ is given by $(1, 0, \delta_1(dx))$ and the Lévy exponent is $\Psi(\xi) = (1 - \cos(\xi)) - i\sin(\xi)$. Let $G = [0, \pi/3]$. Then

$$(\cos(|x|\sin(\xi)))_+ = \cos(|x|\sin(\xi)) \quad \text{and} \quad (\cos(|x|\sin(\xi)))_- \equiv 0,$$

for any $x \in G$. Hence, on the set $G$, we have $g_{1-\beta,-}(x) \equiv 0$ and

$$g_{1-\beta,+}(x) = \int_{\mathbb{R}} e^{-|x|(1-\cos(\xi))} \cos(|x|\sin(\xi))|\xi|^{\beta-1} d\xi$$

$$\geq \tfrac{1}{2} e^{-2|x|} \int_{\mathbb{R}} |\xi|^{\beta-1} d\xi = +\infty.$$

Thus, for any $\mu \in \mathcal{P}(G)$,

$$\int_{\mathbb{R}} \mathcal{E}_{(\Re \chi_\xi)_+}(\mu) \|\xi\|^{\beta-1} d\xi = +\infty \quad \text{and} \quad \int_{\mathbb{R}} \mathcal{E}_{(\Re \chi_\xi)_-}(\mu) \|\xi\|^{\beta-1} d\xi = 0,$$

that is, the conditions of part (e) of the theorem are satisfied, even though $X_t$ is not symmetric since it has a nonvanishing imaginary part of the Lévy exponent.

**Acknowledgment.** The author wishes to thank the anonymous referee for pointing out a serious problem of nonintegrability in the original version of this note.

DEPARTMENT OF MATHEMATICS
UNIVERSITY OF CONNECTICUT
196 AUDITORIUM ROAD
STORRS, CONNECTICUT 06269-3009
USA
E-MAIL: martynas@math.uconn.edu